\documentclass[12pt]{amsart}
\usepackage{amssymb,amsmath,a4wide}

\newcommand{\pa}{\partial}
\newcommand{\ep}{\varepsilon}

\begin{document}

\title[Symmetric-hyperbolic systems]{Symmetric-hyperbolic systems and shock waves}%
\author{Satyanad Kichenassamy}%
\address{Professor of Mathematics and Director,
Laboratoire de Math\'ematiques (UMR 6056), CNRS \&\ Universit\'e
de Reims Champagne-Ardenne, Moulin de la Housse, B.P. 1039,
F-51687
Reims Cedex 2\\France}%
\email{satyanad.kichenassamy@univ-reims.fr}%

\thanks{}%

\date{February 12, 2005}%




\bigskip\medskip

\begin{center}
{\bf SYMMETRIC-HYPERBOLIC SYSTEMS AND SHOCK
WAVES}\footnote{{\bf Encyclopedia of
Mathematical Physics}
(J.-P. Fran\c coise, G. L. Naber, Tsou S. T. Eds.),
Oxford : Elsevier, 2006,  vol. 5, pp. 160-166.}

\bigskip\medskip

{\small SATYANAD KICHENASSAMY}
\end{center}
\bigskip\medskip

\begin{quote}\small{\sc Abstract.}\quad
The theory of symmetric-hyperbolic systems is useful for
constructing smooth solutions of nonlinear wave equations, and for
studying their singularities, including shock waves. We present
the main techniques which are required to apply the theory in
Mathematical Physics.
\end{quote}
\bigskip\medskip

\section{Introduction}

Many systems of partial differential equations arising in
Mathematical Physics and Differential Geometry are
\emph{quasi-linear}: the top-order derivatives enter only
linearly. They may be cast in the form of \emph{first-order}
systems by introducing, if needed, derivatives of the unknowns as
additional unknowns. For such systems, the theory of
\emph{symmetric-hyperbolic} (SH) systems provides a unified
framework for proving the local existence of smooth solutions if
the initial data are smooth. It is also convenient for
constructing numerical schemes, and for studying shock waves.
Despite what the name suggests, the impact of the theory of SH
systems is not limited to hyperbolic problems, two examples being
Tricomi's equation, and equations of Cauchy-Kowalewska type.

Application of the SH framework usually requires a preliminary
reduction to SH form (``symmetrization'').

After comparing briefly the theory of SH systems with other
functional-analytic approaches, we collect basic definitions and
notation. We then present two general rules, for symmetrizing
conservation laws and strictly hyperbolic equations respectively.
We next turn to special features possessed by linear SH systems,
and give a general procedure to prove existence, which covers both
linear and nonlinear systems. We then summarize those results on
shock waves, and on blow-up singularities, which are related to SH
structure. Examples and applications are collected in the last
section.

The advantages of SH theory are: a standardized procedure for
constructing solutions; the availability of standard numerical
schemes; a natural way to prove that the speed of propagation of
support is finite. On the other hand, the symmetrization process
is sometimes \emph{ad hoc}, and does not respect the physical or
geometric nature of the unknowns; to obviate this defect to some
extent, we remark that symmetrizers may be viewed as introducing a
new Riemannian metric on the space of unknowns. The search for a
comprehensive criterion for identifying equations and boundary
conditions compatible with SH structure is still the object of
current research.

The most important fields of application of the theory today are
General Relativity and fluid dynamics, including
magneto\-hydrodynamics.

\subsection{Context of SH theory in modern terms}

The basic reason why the theory works may be summarized as follows
for the modern reader; the history of the subject is, however, more
involved:

Let $H$ be a real Hilbert space. Consider a \emph{linear}
initial-value problem $du/dt+Au=0$; $u(0)=u_0\in H$, where $A$ is
unbounded, with domain $D(A)$. By Stone's theorem, one can solve
it in a generalized sense, if the unbounded operator $A$ satisfies
$A+A^*=0$. This condition contains two ingredients: a symmetry
condition on $A$, and a maximality condition on $D(A)$, which
incorporate boundary conditions (von Neumann, Friedrichs).
Semi-group theory (Hille and Yosida, Phillips, and many others)
handles more general operators $A$: it is possible to solve this
problem in the form $u(t)=S(t)u_0$ for $t>0$, where
$\{S(t)\}_{t\geq 0}$ is a continuous contraction semi-group,
\emph{if and only if} $(Au,u)\geq 0$, and equation $x+Ax=y$ has a
solution for every $y$ in $H$ (this is a maximality condition on
$D(A)$). One then says that $A$ is maximal monotone. For such
operators, $A+A^*\geq 0$. SH systems are systems $Qu_t+Au=F$,
satisfying two algebraic conditions ensuring \emph{formally} that
$A+A^*$ is bounded, and that $Q$ is symmetric and
positive-definite. This algebraic structure enables one to solve
the problem directly, without explicit reference to semi-group
theory. Precise definitions are given next.

We assume throughout that all coefficients, nonlinearities, and
data are smooth unless otherwise specified.

\subsection{Definitions}

Consider a quasi-linear system
\begin{equation}\label{gen}
M_B^{A\alpha}(x,u) \pa_\alpha u^B=N^A(x,u),
\end{equation}
where $u=(u^A)_{A=1,\dots,m}$, $x=(x^\alpha)_{\alpha=0,\dots,n}$,
and $\pa_\alpha=\pa/\pa x^\alpha$. The components of $u$ may be
real or complex. We follow the summation convention on repeated
indices in different positions; $x^0=t$ may be thought of as the
evolution variable; we write $x=(t,\mathbf x)$, with $\mathbf
x=(x^1,\dots,x^n)$. Indices $A$, $B$,\dots\ run from $1$ to $m$,
indices $j$, $k$,\dots\ from $1$ to $n$, and Greek indices from 0
to $n$. The complex conjugate of $u^A$ is written $\bar u^A$.
\begin{itemize}
\item
Equation.~(\ref{gen}) is \emph{symmetrizable} if there are functions
$\sigma_{AB}(x,u)$ such that
\[ M_{AB}^\alpha := \sigma_{AC}M_B^{C\alpha}
\]
satisfies $M_{AB}^\alpha=\bar M_{BA}^\alpha$ for every $\alpha$.
\item
It is \emph{symmetric} if it is symmetrizable with
$\sigma_{AB}=\delta_{AB}$.
\item
It is \emph{symmetric-hyperbolic} with respect to $k_\alpha$ if it
is symmetric and if $k_\alpha M_{AB}^\alpha$ is positive definite:
$k_\alpha M_{AB}^\alpha \xi^A\xi^B>0$ for $\xi=(\xi^A)\neq 0$.
\end{itemize}
Thus, a \emph{symmetrizer} $(\sigma_{AB})$ gives rise to a
Riemannian metric $(k_\alpha \sigma_{AC}M_{B}^{C\alpha})$ on the
space of unknowns, independent of any Riemannian structure on
$x$-space. The system is SH with respect to $x^0$ if
$k_\alpha=\delta_{\alpha}^0$.

The simplest class of SH systems is provided by real
\emph{semi-linear} systems of the form
\begin{equation}
A^0(x)\pa_t u + A^i(x)\pa_i u=N(x,u)
\end{equation}
where the $A^\alpha$ are real symmetric matrices, $A^0$ is
symmetric and positive-definite, and $k_\alpha=\delta_{\alpha 0}$.
Writing $A^0=P^2$, with $P$ symmetric and positive-definite, one
finds that $v=Pu$ solves a SH system with $A^0=I$ (identity
matrix).

\emph{Conservation laws} (with ``reaction'' or ``source'' term
$N^A$) are usually defined as quasi-linear systems of the form
\begin{equation}\label{cl}
\pa_\alpha f^{A\alpha}(x,u)=N^A(x,u).
\end{equation}
They are common in fluid dynamics and combustion. They are
limiting cases of \emph{nonlinear diffusion} equations of the
typical form
\begin{equation} \pa_\alpha
f^{A\alpha}(x,u)=N^A(x,u)+\varepsilon\pa_j(B^{Ajk}_{B}\pa_k u^B).
\end{equation}
The determination of the form of the coefficients $B^{Ajk}_{B}$ is
a non-trivial modeling issue; they may reflect varied physical
processes such as heat conduction, viscosity, or bulk viscosity.
They may depend on $x$, $u$ and the derivatives of $u$. The
simplest case is $B^{Ajk}_{B}=D^{jk}\delta_{B}^A$ with $(D^{jk})$
diagonal. Some authors require the symmetry condition
\begin{equation}\label{bsym}
\delta_{AC}B^{Cjk}_{B}=\delta_{BC}B^{Ckj}_{A}
\end{equation}
Equations in which $f^{A\alpha}=u^A\delta_0^{\alpha}$ are called
\emph{reaction-diffusion} equations; they arise in physical and
biological problems in which chemical reactions and diffusion
phenomena are combined, and in population dynamics.

A conservation law is symmetric if and only if $\pa
f^{A\alpha}/\pa u^B$ is symmetric in $A$ and $B$, which means that
there are, locally, functions $g^\alpha(x,u)$ such that
$f^{A\alpha}=\pa g^{\alpha}/\pa u^A$.

A more fundamental derivation of conservation laws would take us
beyond the scope of this survey.

\section{Symmetrization}

Two general procedures for symmetrization are available: one for
conservation laws, the other for semi-linear strictly hyperbolic
problems.

\subsection{Conservation laws with a convex entropy}

Consider, for simplicity, a conservation law of the form
\begin{equation}\label{cls}
\pa_t u^A+\pa_j f^{Aj}(u)=0.
\end{equation}
We therefore assume that the $f^{A\alpha}=f^{A\alpha}(u)$ and
$f^{A0}(u)=u^A$. We show that the following three statements are
equivalent locally: (i) there is a strictly convex function $U(u)$
such that $\sigma_{AB}=\pa^2 U/\pa u^A\pa u^B$ is a symmetrizer;
(ii) eqn.~(\ref{cls}) implies a scalar relation of the form
$\pa_\alpha U^\alpha=0$, with $U^0$ strictly convex; (iii) there
is a change of unknowns $v_A=v_A(u)$ such that the system
satisfied by $v=(v_A)$ is SH and $(\pa v_A/\pa u^B)$ is
positive-definite.

In fluid dynamics, $U^0$ may sometimes be related to specific
entropy, and $U^j$ to entropy flux. For this reason, if (ii)
holds, one says that $U^0$ is an entropy for  eqn.~(\ref{cls}),
and that $(U^0,U^j)$ is an \emph{entropy pair}. A system may have
several entropies in this sense; this fact is sometimes useful in
studying convergence properties of approximate solutions of
eqn.~(\ref{cls}).

Let us now prove the equivalence of these properties.

Assume first (iii): there are new unknowns $v_A=v_A(u)$ and
functions $g^\alpha(v)$ such that $f^{A\alpha}=\pa g^\alpha/\pa
v_A$. One finds that if eqn.~(\ref{cls}) holds,
\begin{equation}\label{ueq}
\pa_\alpha U^\alpha=0 \text{ where }U^\alpha=v_A\frac{\pa
g^\alpha}{\pa v_A}-g^\alpha.
\end{equation}
Furthermore, we have $f^{A0}=u^A$; therefore, eqn.~(\ref{ueq})
gives: $U^0=v_Au^A-g^0$, so that $U^0$ is the Legendre transform
(familiar from Mechanics) of $g^0$. It follows that $v_A=\pa
U^0/\pa u^A$. Finally, $(\pa v_A/\pa u^B)=(\pa^2 U^0/\pa u^A\pa
u^B)$ is positive-definite, and $U^0$ is strictly convex.

We have proved that (iii) implies (ii).

Next, assume (ii): the \emph{entropy equality} $U_t+\pa_jU^j=0$
holds identically---and not just for the solution at hand. Using
(\ref{cls}), we find
\[ 0=\frac{\pa U}{\pa u^A}\pa_t u^A+\frac{\pa U^j}{\pa u^B}\pa_j u^B
=\left[-\frac{\pa U}{\pa u^A}\frac{\pa f^{Aj}}{\pa u^B}+\frac{\pa
U^j}{\pa u^B}\right]\pa_j u^B.
\]
Assumption (ii) therefore means that $U$ is strictly convex and
satisfies
\begin{equation}\label{assb}
\frac{\pa U}{\pa u^A}\frac{\pa f^{Aj}}{\pa u^B} = \frac{\pa
U^j}{\pa u^B}.
\end{equation}
Now, letting $v_A=\pa U/\pa u^A$ and $g^j(v)=v_Af^{Aj}-U^j$, we
find
\begin{equation}\label{gjeq} \frac{\pa g^j}{\pa v_A}=f^{Aj}+
\left[\frac{\pa U}{\pa u^A}\frac{\pa f^{Aj}}{\pa u^C}-\frac{\pa
U^j}{\pa u^C}\right]\frac{\pa u^C}{\pa v_A}=f^{Aj}.
\end{equation}
Let $\sigma_{AB}={\pa^2 U}/{\pa u^A\pa u^B}$. Since $U$ is
strictly convex, $(\sigma_{AB})$ is positive definite, and so is
its inverse. We have now proved (iii). Note that $u^A=\pa g^0/\pa
v_A$, where $g^0(v)=u^Av_A-U(u)$ is the Legendre transform of $U$.

Next, using eqn.~(\ref{gjeq}), and the relations $\sigma_{AB}=\pa
v_A/\pa u^B=\pa v_B/\pa u^A$, we find
\[ 0=\sigma_{AB}[\pa_t u^B+\pa_jf^{Bj}]
=\sigma_{AB}\pa_t u^B+\frac{\pa v_B}{\pa
u^A}\frac{\pa^2 g^j}{\pa v^B\pa u^C}\pa_j u^C = \sigma_{AB}\pa_t
u^B+\frac{\pa^2 g^j}{\pa u^A\pa u^C}\pa_j u^C,
\]
which is SH; therefore, $\sigma_{AB}$ is a symmetrizer for
eqn.~(\ref{cls}), and (i) is proved.

Thus, (ii) implies (i) and (iii).

Finally, if (i) holds, $\sigma_{AC}\pa f^{Cj}/\pa u^B$ is
symmetric in $A$ and $B$. It follows that
\[ \frac{\pa}{\pa u^C}\left[
 \frac{\pa U}{\pa u^A}\frac{\pa f^{Aj}}{\pa u^B}\right]
 =\sigma_{AC}\frac{\pa f^{Aj}}{\pa u^B}+
 \frac{\pa U}{\pa u^A}\frac{\pa^2 f^{Aj}}{\pa u^B\pa u^C}
\]
is symmetric in $B$ and $C$, so that there are, locally, functions
$U^j$ such that eqn.~(\ref{assb}) holds. Therefore, $(U,U^j)$ is
an entropy pair, and we see that (i) implies (ii).

This completes the proof of the equivalence of (i), (ii) and
(iii).

\subsection{Strictly hyperbolic equations}

Consider the scalar equation $Pf=g(t,\mathbf x)$, where $P$ is the
linear operator $P=\pa_t^N-\sum_{j=0}^{N-1} p_{N-j}(t,\mathbf
x)\pa_t^j$, of order $N$. Let $\Lambda=(1-\Delta)^{1/2}$, where
$\Delta$ is the Laplace operator on the space variables. Then
$u=(u^A)$, where $u^A=\pa_t^{A-1}\Lambda^{N-A}f$ for $A=1,\dots,
N$, solves a first-order \emph{pseudo-differential} system of the
form
\[ u_t-Lu=G.
\]
If $P$ is strictly hyperbolic, the principal symbol $a_1(t,\mathbf
x,\mathbf\xi)$ of $L$ has a diagonal form with real eigenvalues
$\lambda_j(t,\mathbf x,\mathbf\xi)$, and there are projectors
$p_j(t,\mathbf x,\mathbf\xi)$ ($p_j^2=p_j$) which commute with
$a_1$, such that $1=\sum_j p_j$, and $a_1=\sum_j \lambda_j p_j$.
Let $r_0=\sum_j p_j^*p_j$, and $r_0(D)$ the corresponding
operator. Equation
\[ r_0(D) \pa_tu -r_0(D)Lu=r_0(D)G
\]
is formally SH in the following sense: (i) $r_0$ is
positive-definite and $r_0 a_1$ is Hermitian.

\section{Linear problems}

Consider a linear system
\begin{equation}
Lu:=Q(t,\mathbf x)\pa_t u
+A^j(t,\mathbf x)\pa_j u + B(t,\mathbf x)u=f(t,\mathbf x).
\end{equation}
We assume that $Q$ and the $A^j$ are real and symmetric, $Q\geq c$
with $c$ positive, and all coefficients and their first-order
derivatives are bounded.

\subsection{Energy identity}
Multiplying the equation by $u^T$ (transpose of $u$), one derives
the ``energy identity''
\begin{equation}\label{energy}
\pa_t (u^TQu)+\pa_j (u^TA^ju)+u^TCu=2u^Tf(t,\mathbf x),
\end{equation}
where $C=2B-\pa_tQ-\pa_jA^j$. $C$ is not necessarily positive.
However, $v:=u\exp(-\lambda t)$ satisfies a linear SH system for
which $C$ is positive-definite if $\lambda$ is large enough.

\subsection{Propagation of support}

A basic property of wave-like equations is finite speed of
propagation of support: if the right-hand side vanishes, and if
the solution at time 0 is localized in the ball of radius $r$,
then the solution at time $t$ is localized in the ball of radius
$r+ct$ for a suitable constant $c$.

This property also holds for SH systems. To see this, let us
consider the set where a solution $u$ vanishes: if the initial
condition vanishes for $|\mathbf x|\leq R$, we claim that $u$ at
some later time vanishes for $|\mathbf x|\leq R-t/a$, for $a$
large enough.

Indeed, let us integrate the energy identity on a truncated cone
$\Gamma :=\{ |\mathbf x|\leq a(t_0-t)/t_0; 0\leq t\leq t_1\}$ with
$t_1<t_0$. The boundary of $\Gamma$ consists of three parts:
$\pa\Gamma = \Omega_0\cup\Omega_1\cup S$, where $\Omega_0$ and
$\Omega_1$ represent the portions of the boundary on which $t=0$
and $t_1$ respectively. The outer normal to $S$ is proportional to
$(a,t_0x^j/|\mathbf x|)$. Let $E(s)$ denote the integral of
$u^TQu$ on $\Gamma\cap\{t=s\}$. Integrating eqn.~(\ref{energy}) by
parts, we obtain
\begin{equation}
E(t_1)-E(0)+\int_S u^T\Phi u\, ds
=\iint_\Gamma(2u^Tf-u^TCu)dt\,d\mathbf x,
\end{equation}
where $\Phi$ is proportional to $aQ+t_0\sum_jx^jA^j/|\mathbf x|$.
Take $a$ so large that $\Phi$ is positive-definite. The integral
over $S$ is then non-negative. If $C$ is positive-definite and
$f\equiv 0$, so that $E(0)=0$, we find that $E(t_1)\leq 0$. Since
$Q$ is positive-definite, this implies $u\equiv 0$ on $\Omega_1$,
as claimed.

\subsection{A numerical scheme}
System $Lu=f$ may be discretized, \emph{e.g.} by the
Lax-Friedrichs method: Let $h$ be the discretization step in
space, and $k$ the time step; write $\tau_ju(t,\mathbf
x)=u(t,x^1,\dots,x^j+h,\dots,x^n)$ (translation in the $j$
direction). one replaces $\pa_j u$ by the centered difference in
the $j$ direction: $(\tau_ju-\tau_j^{-1}u)/2h$; and the time
derivative by
\begin{equation}
[u(t+k,\mathbf x)-
\frac1{2n}\sum_j(\tau_j u(t,\mathbf x)+\tau_j^{-1}u(t,\mathbf x))]/k.
\end{equation}
For consistency of the scheme, we require $k/h=\lambda>0$ to be
fixed as $k$ and $h$ tend to zero; stability then holds if
$\lambda$ is small.

\section{Nonlinear problems and singularities}

We give a simple set-up for proving the existence of smooth
solutions to SH systems for small times. Such solutions may
develop singularities. We limit ourselves to two types of
singularities, on which SH structure provides some information:
jump discontinuities and blow-up patterns. Caustic formation is
not considered.

\subsection{Construction of a smooth solution}

Consider a real SH system (eqn.~(\ref{gen})). Recall that a
function of $\mathbf x$  belongs to the Sobolev space $H^s$ if its
derivatives of order $s$ or less are square-integrable. One
constructs a solution defined for $t$ small, which is in $H^s$,
$s>n/2+1$ as a function of $\mathbf x$, by the following
procedure:
\begin{itemize}
\item[(1)]
Replace spatial derivatives by regularized
operators, which should be bounded in Sobolev spaces; the
regularized equation is an ODE in $H^s$; let $u_\ep$ be its
solution.
\item[(2)]
Write the equation satisfied by derivatives of
order $s$ of $u_\ep$, and apply the energy identity to it.
\item[(3)]
Find a positive $T$ such that the solution is bounded in $H^s$ for
$|t|\leq T$, uniformly in $\ep$; this implies a $C^1$ bound.
\item[(4)]
Prove the convergence of the approximations in $L^2$.
\item[(5)]
Prove the continuity in time of the $H^s$ norm; conclude that the
$u_\ep$ tend to a solution in $C(-T,T;H^{s})$.
\end{itemize}
The result admits a local version, in which Sobolev spaces are
replaced by Kato's ``uniformly local'' spaces. Uniqueness of the
solution is proved along similar lines. We do not attempt at
identifying the infimum of the values of $s$ for which the Cauchy
problem is well-posed.

\subsection{Jump discontinuities: Shock waves}

A \emph{shock wave} is a weak solution of a system of conservation
laws admitting a jump discontinuity. By definition, weak solutions
satisfy, for any smooth function $\phi_A(x)$ with compact support,
$\iint \{f^{A\alpha}\pa_\alpha\phi_A+N^A\phi_A\} dt\, d\mathbf
x=0$.

The theory of shock waves is an attempt to understand solutions of
conservation laws which are limits of solutions of diffusion
equations; the hope is that the influence of second-derivative
terms is appreciable only near shocks, and that, for given initial
data, there is a unique weak solution of the conservation law
which may be obtained as such a limit, if modeling has been done
correctly. This problem may be difficult already for a single
shock (``shock structure'').

The theory of shock waves follows the one-dimensional theory
closely. We therefore describe the main facts for a conservation
law in one space dimension ($u=u(t,x)$):
\[ \pa_tu+\pa_xf(u)=0.
\]
If a shock travels at speed $c$, the weak formulation of the
equations gives the \emph{Rankine-Hugoniot relation}
$c[u]=[f(u)]$, where square brackets denote jumps. There may be
several weak solutions having the same initial condition. One
restricts solutions by making two further requirements: (a) the
system admits is an entropy pair $(U,F)$ with a convex entropy and
(b) to be admissible, weak solutions must be limits of ``viscous
approximations''
\[ \pa_tu+\pa_xf(u)=\ep \pa^2_{x}u
\]
as $\ep\to 0$. One then finds easily that the entropy equality
($\pa_t U+\pa_x F=0$) must be replaced, for such weak solutions,
by the \emph{entropy condition}: $\pa_t U+\pa_x F\leq 0$ in the
weak sense. This condition admits a concrete interpretation if the
gradient of each characteristic speed is never orthogonal to the
corresponding right eigenvector (``genuine nonlinearity''); in
that case, characteristics must impinge on the shock (``shock
inequalities'').

For the equations of gas dynamics with polytropic law
($pv^\gamma=$ const.), there is a unique solution with initial
condition $u=u_l$ for $x<0$, $u=u_r$ for $x>0$, where $u_l$ and
$u_r$ are constant (``Riemann problem'') which satisfies the
entropy condition, provided $|u_l-u_r|$ is small. More generally,
if the equation of state $p=p(v,s)>0$ satisfies $\pa p/\pa v<0$
and $\pa^2 p/\pa v^2>0$, the shock inequalities are equivalent to
the fact that the entropy increases after the passage of a shock
with $|u_l-u_r|$ small.

On the numerical side, one should mention: (i) the widely used
idea of upstream differencing; (ii) the Lax-Wendroff scheme, the
complete analysis of which requires tools from soliton theory;
(iii) the availability of general results for dissipative schemes
for SH systems.

Recent trends include: (i) admissibility conditions when genuine
nonlinearity does not hold; (ii) other approximations of shock
wave problems, most notably kinetic formulations.

Some of the ideas of shock wave theory have been applied to
Hamilton-Jacobi equations and to motion by mean curvature, with
applications to front propagation problems and Computer Vision.

\subsection{Stronger singularities: Blow-up patterns}

The amplitude of a solution may also grow without bound. Examples
include optical pulse propagation in Kerr media and singularities
in General Relativity. The phenomenon is common when reaction
terms are allowed. As we now explain, this phenomenon is reducible
to SH theory in many cases of interest.

Blow-up singularities are usually not governed by the
characteristic speeds defined by the principal part, because
top-order derivatives are balanced by lower-order terms. In many
applications, a systematic process (Fuchsian Reduction) enables
one to identify the correct model near blow-up; as a result, one
can write the solution as the sum of a singular part, known in
closed form, and a regular part. If the singularity locus is
represented by $t=0$, the regular part solves a renormalized
equation of the typical form
\begin{equation}\label{fs}
tMu+Au=t^\ep N,
\end{equation}
where $Mu=0$ is SH. Under natural conditions, for any initial
condition $u_0$ such that $Au_0=0$, there is a unique solution of
eqn.~(\ref{fs}) defined for small $t$.

The upshot is an asymptotic representation of solutions which
renders the same services as an exact solution, and is valid
precisely where numerical computation breaks down.

Fuchsian Reduction enables one in particular to study (i) the
blow-up time; (ii) how the singularity locus varies when Cauchy
data, prescribed in the smooth region, are varied; (iii)
expressions which remain finite at blow-up. It is the only known
general procedure for constructing analytically singular
space-times involving arbitrary functions, rather than arbitrary
parameters, and is therefore relevant to the search for
alternatives to the big-bang.

\section{Examples and applications}

\subsection{Wave equation with variable coefficients}

Consider equation
\[ \pa_{tt}u+2a^j(x)\pa_{jt}u-a^{jk}(x)\pa_{jk}u=f(t,\mathbf x,u,\nabla u),
\]
with $(a^{jk})$ positive-definite. Letting
$v=(v_0,\dots,v_{n+1}):=(u,\pa_j u,\pa_t u)$, we find the system
\begin{align*}
\pa_t v_0 &= v_{n+1},\\ \pa_t v_k - \pa_k v_{n+1} &=0,\\ \pa_t
v_{n+1} +2a^k\pa_k v_{n+1}-a^{jk}\pa_{k}v_j &= f.
\end{align*}
It is symmetrizable, using the quadratic form
$\sigma_{AB}v^Av^B=v_0^2+a^{jk}v_jv_k + v_{n+1}^2$.

One proves directly that if $v_j=\pa_j v_0$ for $t=0$, this
relation remains true for all $t$.

\subsection{Maxwell's equations}

Maxwell's equations may be split into six evolution equations:
$\pa_t \mathbf E-\mathop{{\rm curl}} \mathbf B+\mathbf{j}=0$ and
$\pa_t \mathbf B+\mathop{{\rm curl}} \mathbf E=0$, and two
``constraints'' $\mathop{{\rm div}} \mathbf E-\rho=0$,
$\mathop{{\rm div}} \mathbf B=0$. The system of evolution
equations is already in symmetric form; the quadratic form
$\sigma_{AB}u^Au^B$ is here $|\mathbf E|^2+|\mathbf B|^2$.

\subsection{Compressible fluids}

Consider first the case of a polytropic gas:
\begin{equation}
\begin{array}{rl}
\pa_t\mathbf v+(\mathbf v\cdot\nabla)\mathbf v+\rho^{-1}\nabla
p&=0\quad \\
\pa_t\rho+\mathop{{\rm div}}(\rho\mathbf v)&=0
\end{array}
\end{equation}
with $p$ proportional to $\rho^\gamma$. Taking $(p,\mathbf v)$ as
unknowns, one readily finds the SH system
\begin{gather}
\frac1{\gamma p}\pa_t p+\frac1{\gamma p}(\mathbf v\cdot\nabla)
p+\mathop{{\rm div}}\mathbf v=0\\ \rho\pa_t \mathbf v+\nabla
p+\rho (\mathbf v\cdot\nabla) \mathbf v=0.
\end{gather}

Symmetrization for more general compressible fluids with
dissipation, including bulk viscosity, so as to satisfy the
additional condition (\ref{bsym}) may be achieved if we take as
thermodynamic variables $\rho$ and $T$, and assume pressure $p$ and
internal energy $\ep$ satisfy $\pa p/\pa\rho>0$ and $\pa \ep/\pa
T>0$, by taking as unknowns $(\rho,\rho \mathbf v,
\rho(\ep+|\mathbf v|^2/2))$. The specific entropy $s$ satisfies
$d\ep=Tds-pd(1/\rho)$. If the viscosity and heat conduction
coefficients are positive, one finds that $U=-\rho s$ is a convex
entropy (in the sense of SH theory) on the set where $\rho>0$,
$T>0$.

\subsection{Einstein's equations}

The computation of solutions of Einstein's equations over long
times, in particular in the study of coalescence of binary stars,
has recently led to unexplained difficulties in the standard
Arnowitt-Deser-Misner (ADM) formulation of the initial-value
problem in General Relativity. One way to tackle these
difficulties is to re-write the field equations in SH form; we
focus on this particular aspect of recent research.

Recall the problem: find a four-dimensional metric $g_{ab}$ with
Lorentzian signature, such that $R_{ab}-\frac12 Rg_{ab}=\chi
T_{ab}$, with $\nabla^aT_{ab}=0$, combined with an equation of
state if necessary. $R_{ab}$ is the Ricci tensor and
$R=g^{ab}R_{ab}$ is the scalar curvature; they depend on
derivatives of the metric up to order two. In addition to the
metric, $T_{ab}$ involves physical quantities such as fluid
four-velocity or an electromagnetic field. The conservation laws
of classical Mathematical Physics are all contained in the
relation $\nabla^aT_{ab}=0$.

Now, the field equations cannot be solved for $\pa^2_tg_{ab}$,
and, as a consequence, the Taylor series of $g_{ab}$ with respect
to time cannot be determined, even formally, from the values of
$g_{ab}$ and $\pa_tg_{ab}$ for $t=0$ (\emph{i.e.}, the Cauchy
data). Furthermore, these data must satisfy four \emph{constraint
equations}. If the constraints are satisfied initially, they
``propagate.'' But in numerical computation, these constraints are
never exactly satisfied, and the computed solution may deviate
considerably from the exact solution. Also, numerical computations
depend heavily on the way Einstein's equations are formulated.

The simplest way to derive a SH system is to replace $R_{ab}$ by
$R^{(h)}_{ab}=R_{ab}-\frac12[g_{bc}\pa_aF^c+g_{ac}\pa_bF^c]$,
where $F^c:=g^{ab}\Gamma^c_{ab}$. It turns out that
$R^{(h)}_{ab}=-\frac12 g^{cd}\pa_{cd}g_{ab}+H_{ab}(g,\pa g)$,
where the expression of $H_{ab}$ is immaterial. Applying to each
component of the metric the treatment of the first example above
(wave equation with variable coefficients), one easily derives a
SH system of 50 equations for 50 unknowns: the 10 independent
components of the metric, and their 40 first-order derivatives.
Now, if the $\Gamma^c$ are initially zero (coordinates are
``harmonic''), they remain so at later times.

Unfortunately, the harmonic coordinate condition does
not seem to be stable in the large. More recent formulations start
with one of the standard set-ups (ADM formalism, conformal
equations, tetrad formalism, Newman-Penrose formalism) and proceed
by adding combinations of the constraints to the equations,
multiplied by parameters adjusted so as to ensure hyperbolicity or
symmetric-hyperbolicity if needed. Another recent idea is to add a
new unknown $\lambda$ which monitors the failure of the constraint
equations; one adds to the equations a new relation of the form
$\pa_t\lambda = \alpha C-\beta\lambda$, where $C=0$ is equivalent
to the constraints, and $\alpha$ and $\beta$ are parameters. One
then adds coupling terms to make the extended system SH. It is
expected that the set of constraints acts as an \emph{attractor}.

Reported computations indicate that these methods have resulted in
an improvement of the time over which numerical computations are
valid.

\subsection{Tricomi's equation}

Let $\varphi(x,y)$ solve $(y\pa_x^2-\pa_y^2)\varphi=0$. Letting
$u=e^{-\lambda x}(\pa_x\varphi,\pa_y\varphi)$, one finds a
symmetric system $Lu=0$, with
\[ L=\begin{pmatrix}y & 0 \\ 0 & 1\end{pmatrix}(\pa_x+\lambda)
-\begin{pmatrix}0 & 1 \\ 1 & 0\end{pmatrix}\pa_y.
\]
If $Z= \begin{pmatrix}1 & y \\ 1 & 1\end{pmatrix}$, we find that
$K=ZL=A^1\pa_x+A^2\pa_y+B$, where $B-\frac12(\pa_x A^1+\pa_y A^2)=
\begin{pmatrix}1/2+\lambda y & \lambda y \\
                 \lambda y & \lambda\end{pmatrix}$
is positive definite if $y$ is bounded, of arbitrary sign, and
$\lambda$ is small.

\subsection{Cauchy-Kowalewska systems}

Consider a complex system
\begin{equation}\label{ck}
\pa_t u = A^j(z,t,u)\frac{\pa u}{\pa z^j}+B(z,t,u),
\end{equation}
where $u=(u^A)$, $z=(z^1,\dots, z^n)$. The coefficients are
analytic in their arguments when $z$ and $t$ are close to the
origin and $u$ is bounded by some constant $K$. The
Cauchy-Kowalewska theorem ensures that, for any analytic initial
condition near the origin, this system has a unique analytic
solution near $z=0$, \emph{even without any symmetry assumption on
the} $A^j$. This result is a consequence of SH theory
(Garabedian).

Indeed, write $z^j=x^j+iy^j$,
$\pa_{z^j}=\frac12(\pa_{x^j}-i\pa_{y^j})$, and $\pa_{\bar
z_j}=\frac12(\pa_{x^j}+i\pa_{y^j})$. Recall that analytic
functions of $z$ satisfy the Cauchy-Riemann equations $\pa_{\bar
z_j}u=0$.

Adding $(\bar A^j)^T\pa_{\bar z_j}$ to (\ref{ck}), and using the
definition of $\pa_{z^j}$ and $\pa_{\bar z_j}$, we find the
symmetric system
\begin{equation}\label{cks} u_t = \frac12(A^j+(\bar A^j)^T)\pa_{x^j}u
+\frac1{2i}(A^j-(\bar A^j)^T)\pa_{y^j}u+B.
\end{equation}
Solving this system, we find a candidate $u$ for a solution of
eqn.~(\ref{ck}). To show that $u$ is analytic if the data are, we
solve a second SH system for $w=w^{(j)}:=\pa_{\bar z_j}u$. If the
data are analytic, $w$ vanishes initially, and therefore remains
zero for all $t$. Therefore, $u$ is indeed analytic.


\section*{See also}

Linear and nonlinear evolution equations. Computational methods in
General Relativity. Initial-value problem for Einstein's
equations. Magneto\-hydrodynamics. PDE in fluid mechanics.

\section*{Keywords}

Nonlinear wave equations, symmetrizer, hyperbolic, conservation
laws, Cauchy problem, singularities, entropy, hydrodynamics,
Relativity, Tricomi problem.

\end{document}